\documentclass[3p,authoryear,final]{elsarticle}    
\usepackage{amsmath,amsfonts,amssymb}
\usepackage{graphicx,psfrag,epsf}
\usepackage{natbib}
\usepackage{tikz}
\usepackage{pgfplots}
\usepackage[T1]{fontenc}
\usepackage{multirow}
\usepackage{color}
\usepackage{soul}

\newcommand{\RR}{\mathbb{R}}

\newcommand{\CC}{{\cal C}}

\newcommand{\pentiere}[1]{\lfloor #1 \rfloor}

\newcommand{\pto}[2]{\langle #1,#2 \rangle}

\date{}
\begin{document}
\begin{frontmatter}
\title{SOCP relaxation bounds for the optimal subset selection problem  applied to robust linear regression \tnoteref{t1}}
\tnotetext[t1]{This manuscript has been accepted for publication in European Journal of Operational Research. The manuscript will undergo 
copyediting, typesetting, and review of the resulting proof before it is published in its final form. 
Please note that during the production process errors may be discovered which could affect the content, and all disclaimers that apply to 
the journal apply to this manuscript. A definitive version was subsequently published in \\ European Journal of Operational Research {\bf 246} 
(2015) pp. 44-50. DOI:10.1016/j.ejor.2015.04.024}

\author[yo]{Salvador Flores}
\fntext[f1]{Supported by Fondecyt under grant 3120166.}
\address[yo]{Centro de Modelamiento Matem\'atico (CNRS UMI 2807)--Universidad de Chile. \\
Beauchef 851, Santiago, Chile.}
\ead{sflores@dim.uchile.cl}

\begin{abstract}
{ This paper deals with the problem of finding the globally optimal subset of $h$ elements 
from a larger set of $n$ elements in $d$ space dimensions so as to minimize a quadratic criterion, with an special emphasis on 
applications to computing  the Least Trimmed Squares Estimator (LTSE) for robust regression.}
The computation of the LTSE is a challenging  subset selection problem involving a 
nonlinear program with continuous and binary variables, linked in a highly nonlinear fashion.
The selection of a globally optimal subset using the branch and bound (BB) algorithm is 
limited to problems {in very low dimension, tipically $d\leq 5$, as the complexity of the problem increases 
exponentially with $d$.} We introduce a bold pruning strategy in the BB algorithm that results in a significant reduction 
in computing time, at the price of a negligeable accuracy lost.
The novelty of our algorithm is that the bounds at nodes of the BB tree come from pseudo-convexifications derived using
a linearization technique with approximate bounds for the nonlinear terms. The approximate bounds are computed solving an
auxiliary {semidefinite optimization} problem. 
We show through a computational study that our algorithm performs well in a wide set of
the most difficult instances of {the LTSE}  problem. 

\end{abstract}

\begin{keyword}
Global Optimization; Integer Programming; High Breakdown Point Regression; Branch and bound; Relaxation-Linearization Technique.
\end{keyword}
\end{frontmatter}

\section{Introduction}
In this article we deal with a nonlinear subset selection problem arising 
in the  computation of linear regression estimators with strong robustness properties.

\begin{figure}
\includegraphics[scale=.8]{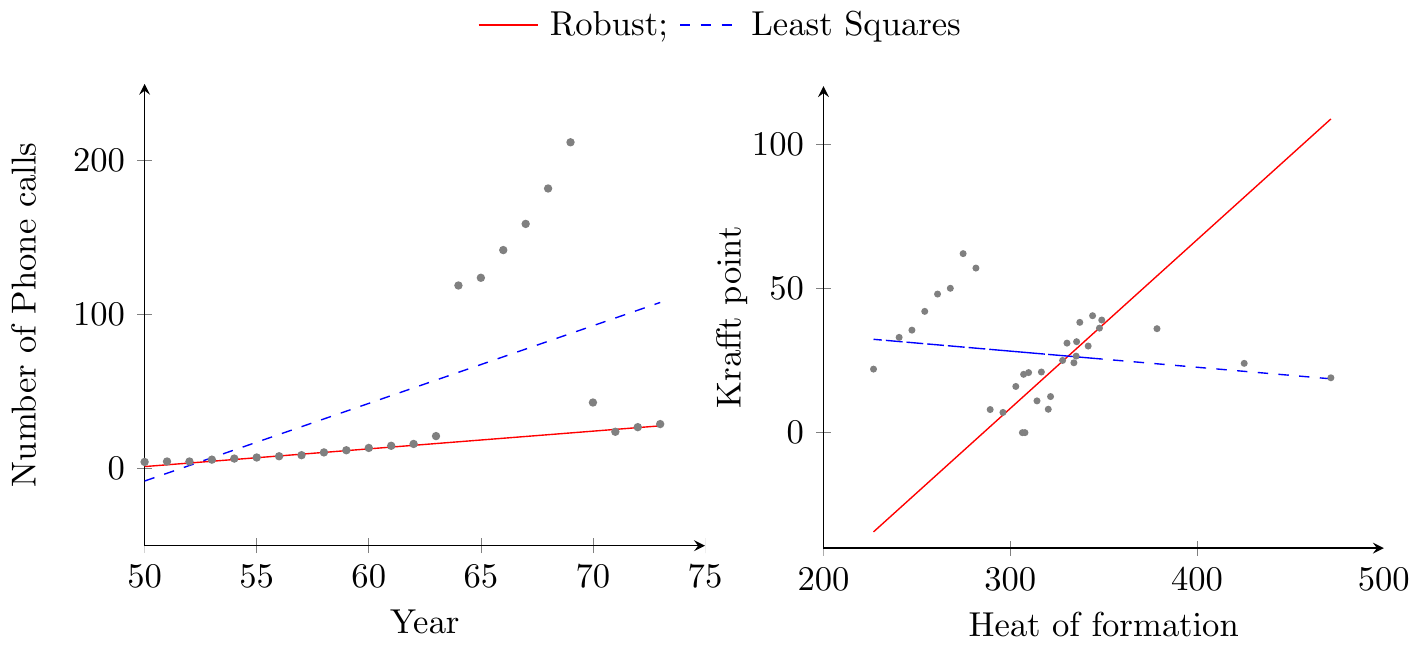}
 \caption{Some real data containing outliers. At the left: phone calls from Belgium; at the right: Krafft point of chemical compounds.}\label{F:Datasets}
\end{figure}

Before entering into the details of the problem we introduce the problem of robust estimation through an example from \citet{Rousseeuw-Leroy87}.
In Figure \ref{F:Datasets} (left) we show, for each year from 1950 to 1973, the number of outgoing international phone calls from Belgium. 
The bulk of the data follows a linear model; nonetheless, there are 6 observations that deviate from  the majority. 
In fact, during the period between 1964 and 1969, there was a change on the record system, which actually 
recorded the total duration, in minutes, of the international phone calls instead of the number of calls.
We plot the regression line obtained by a robust method (solid line) and that obtained by the 
method of Least Squares (dashed line). The Least Squares estimation is strongly affected by the outliers.
 
\pgfplotsset{every axis legend/.append style={
        at={(0.01,1)}, 
        anchor=north west}}

A more appealing example \citep{krafft} where a subpopulation acts differently, is shown in 
Figure \ref{F:Datasets} (right). 
There are plotted, for 32 chemical compounds, a quantity called  Krafft point versus a molecular descriptor 
called heat of formation. There is a main group that follows a regression line, correctly estimated by a robust 
estimator (solid line); there is also, besides some few outliers at right, a second, smaller group forming what seems 
to be another regression line. The observations in the second group correspond to sulfonates. The 
Least Squares estimator (dashed line) is not helpful to detect the presence of the second group.

As the reader can figure out, in higher dimensions, where visual inspection is not longer 
an alternative for detecting outliers, specifically suited methods are needed to deal with outliers.
This is what robust estimators are about.
Besides robustness, in a sense to be specified soon, robust regression estimators satisfy statistical properties 
such as asymptotic normality, square-root rate of convergence and equivariance properties \citep{Maronna-Martin-Yohai2006}. 
Unfortunately, the use of robust estimators is not as widespread as one may expect because their computation 
is very time-consuming. Unlike other problems arising in Statistics, the difficult problems involved 
in computing robust estimators remain almost unknown to O.R. specialists.

\subsection{Description of the {motivating} problem}

We have at our disposal a sample consisting of $n$  observations of the explicative variables $\{x_1,...,x_n\} \subseteq \RR^{d}$.
To each observation of explicative variables it corresponds a response or dependent variable, gathered together in a vector $y=(y_1,...,y_n)\in \RR^n$.
We assume that the random variables $x$ and $y$ are related through a linear model, which implies the existence of a 
vector $\beta\in \RR^d$ such that
\begin{equation}\label{E:linmod}
 y_i = x_{i}^\top \beta +\delta_i,
\end{equation}

with $\delta_i$ i.i.d. , $E[\delta_i ]=0$, $Var[\delta_i]=\sigma^2$. 
The objective of linear regression is to estimate the parameter $\beta$.

For convenience, we put the explicative variables as rows of a matrix $X\in \RR^{n\times d}$,
\begin{equation*}
X=\left( \begin{matrix}
x_1\\ x_2 \\ \vdots \\ \vdots \\ x_n
         \end{matrix}
    \right) =\left( \begin{matrix}
x_1^{(1)}& x_1^{(2)}& \cdots&x_1^{(d)}\\
x_2^{(1)}& x_2^{(2)}& \cdots&x_2^{(d)} \\ \vdots&\vdots&\vdots&\vdots \\ \vdots&\vdots&\vdots&\vdots \\
x_n^{(1)}& x_n^{(2)}& \cdots&x_n^{(d)}
         \end{matrix}
    \right).
 \end{equation*}
For $\beta\in \RR^d$, we denote by $\mathbf r(\beta)$ the vector of \emph{residuals} 
$\mathbf r=y-X\beta$, with components $r_i=y_i- x_i^\top\beta$. 
The Least Squares (LS) estimator is obtained by minimizing the sum of the squared residuals:
\begin{equation*}
\min\limits_{\beta\in\RR^d} \sum_{i=1}^n r_i(\beta)^2.  
\end{equation*}

Hereafter we adopt the robustness notion introduced by \citet{Donoho-Huber83}, which is based on the concept of \emph{Breakdown Point (BDP)}.
The BDP of an estimator on a sample is defined as the minimum fraction of  observations that
need to be replaced by arbitrary ones for the estimator to take on arbitrary values. The BDP of the common LS estimator is $1/n$, since it suffices 
to control just one observation to make the estimator divergent.
Therefore, the asymptotic BDP of the LS estimator as $n$ grows to infinity is $0\%$. The same is true if the LS estimator is 
replaced by any estimator obtained by minimizing a convex function of the residuals.
Since the pioneer Least Median of Squares (LMS) estimator \citep[see][for a precise description]{Rousseeuw-Leroy87}, there has been a continuous 
improvement leading to high BDP estimators with optimal statistical properties, such as asymptotic normality, speed of convergence and efficiency. 
In this article we focus on the Least Trimmed Squares estimator, which has the best statistical properties and is defined through a well structured 
optimization problem.
Let $h$ be an integer number comprised between $n/2$ and $n$, and denote by $|r|_{1:n}\leq |r|_{2:n} \leq \cdots \leq |r|_{n:n}$  the residuals, 
ordered by increasing absolute value. 
 The Least Trimmed Squares (LTS) estimator is defined as a solution of the problem:
\begin{equation*}\label{P:LTS} 
\min_{\beta \in \RR^d} \sum_{i=1}^h  r(\beta)_{i:n}^2. \tag{LTS}
\end{equation*}
In words, the LTSE is the vector of regression coefficients $\hat \beta$ that minimizes the sum of the $h$ smallest 
squared distances from the hyperplane defined by $\hat\beta$ to the observations $y_i$. 
The LTSE attains the maximum asymptotic BDP of $50\%$ 
by taking $h=\pentiere{n/2}+\pentiere{(d+1)/2}$ \citep{Rousseeuw-Leroy87}.

\subsection{State of the art}

The first approaches to the optimal subset selection problem appeared in the field of pattern recongnition 
\citep{Narendra-Fukunaga,YUAN1993,CHEN2003,Somol2004}. Unlike robust regression, the feature selection problem 
addressed there is a maximization problem, whose difficulties are somehow different from ours.
To the best of our knowledge, the LTSE is the only reported application of subset selection involving minimization.

 The { exact} computation of high-BDP estimators { for $d$ greater than, say, $5$} is a difficult global optimization problem
\citep{Thorsten,Mount2006,Mount2014}.
{
Indeed, \citet{Mount2006} provide results suggesting that any exact algorithm for the related LMS requires, for large $n$, 
a time superior to $k\cdot n^d$ for some constant $k>0$. \citet{Mount2014} extend this result to the LTSE under 
consideration here  proving that, up to a constant, the time required for computing the LTSE for a given 
coverage level $h$ must be, for large $n$, bounded between $(n/h)^d$ and $n^{d+1}$.  
}
 
For this reason, the overwhelming majority of the literature on computation of robust estimators is focused on stochastic approximation algorithms. 
Most of these algorithms are constructed upon the basic resampling algorithm proposed originally by \citet{Rousseeuw-Leroy87} for computing the LMS
estimator. \citet{Rousseeuw-VanDriessen2006} developed a refined version including local improvements and adapted to the LTSE. 
Recently, \citet{TortiEtAl2012} conducted a benchmark of stochastic algorithms for {approximating} high-BDP linear regression estimators. 
The interested reader can find in that article an up-to-date account of stochastic approximation algorithms for 
{robust (though not high-BDP) linear} regression.
{
Stochastic approximation algorithms tipically provide good approximations to the actual estimator for problems 
with a number of observations in the hundreds or even in the few thousands. However, they have some shortcomings as well.
Two different runs or different implementations of the algorithm may give different results. 
Also, as it is not guaranteed that the obtained solution is a global minimum, nothing can 
be said about the actual breakdown point of such approximations. Even a deterministic, 
constant-factor approximation to a high breakdown point estimator may have a very low breakdown point.

That being said, for small or medium-sized datasets one may be disposed to invest more time to have 
a guaranteed global solution in return. Unfortunately, despite the notable literature dedicated to stochastic approximations there exist very few 
deterministic algorithms for computing robust estimators, exactly or approximately.} We can mention the proposals by \citet{Steele-Steiger1986} and 
\citet{Stromberg1993} for computing the LMS estimator; both based on enumeration of elemental subsets. { For the particular case of LMS 
regression with two predictors ($d=2$), \citet{Mount2007} devised a Branch and Bound algorithm with an asymptotic running time of $O(n \log_2 n)$.}

The first great step forward {in the computation of the LTSE} came with the Branch and Bound algorithm (BBA) of \citet{Agullo2001}.
The BBA is an adaptation {for minimization problems} of the  \textquoteleft feature subset selection' branch and bound maximization 
algorithm by \citet{Narendra-Fukunaga} and relies on the monotonicity of the problem.
At a glance, the BBA enumerates all the subsets of $h$ observations out of the 
$n$, by starting from the empty set and adding one observation at a time. Since the sum of the squared residuals increases when an additional 
observation is added to the LS fit, if a subset of observations with cardinality $k<h$ is found to give a sum of squared 
residuals larger than that of the incumbent set, then by monotonicity all the sets containing that set can be discarded from further examination. 
The BBA is reported to be efficient in datasets with up to about $n=50$ observations and $d=5$ features. 
To a large extent, the difficulty of the BBA in tackling larger problems stems from the following limitations of the monotonicity bound
\begin{itemize}
\item[\textbf{--}]
{ It does not provide useful information for small subsets:} since in dimension $d$ it is always possible to fit $d$ observations exactly,
it is impossible to obtain a non-trivial (positive) lower bound for the sum of squared residuals for a fit of $h$ observations containing a given 
subset of $d$ or less observations. In the enumeration tree, this amounts to not having a lower bound to prune at the top $d$  levels.  
\item[\textbf{--}]
{It does not look ahead:} for instance, if $n = 15$ and $h = 6$ it gives ``the sum of the squared
residuals of a regression over six observations comprising observations 2, 5 and 8 is greater than the sum of the squared residuals of the regression
over observations 2, 5 and 8'' without quantifying the increase in sum of squared residuals due to the incorporation of three
more observations. 
\end{itemize}

Consequently, as the dimension of the explicative variables increases, the BBA must examine a large number of elemental subsets, since pruning
is possible only at the bottom levels of the enumeration tree, even if a good global upper bound is available beforehand.
On the other hand, the BBA uses a quite efficient explicit formula for computing the increase in sum of squared residuals when adding one 
observation, which makes his algorithm quite efficient at the lower levels of the enumeration tree.

\citet{Gatu2010} extended the BBA for obtaining the LTSE for many coverage values $h$ at once, besides 
 improving the numerical linear algebra used to compute lower bounds.
{Very recently, \citet{Bertsimas} proposed a linear Mixed Integer Optimization (MIO) formulation of the Least Quantile of Squares
(and in particular the LMS) regression problem. They report good results at solving to provable optimality problems of small ($n=100$) and medium 
($n=500$) size for $d=3$.
Some impressive results are reported for approximate (albeit deterministic) solutions for problems with $d\leq 20$ and $n$ in the order of $10^4$. 
The BDP of the regression performed with approximate solutions is not reported.
}

\subsection{Innovations and contributions}
We model the computation of the LTSE as a nonconvex optimization problem comprising continuous and binary variables with nonlinear 
interdependence. Since the nonlinear coupling of the variables makes obtaining bounds on the continuous variables impractical, we devise a technique 
to obtain approximate bounds on the continuous variables; this is done by solving an Semi-Definite Programming (SDP) problem only once. 
Then we use the approximate bounds to obtain, via { the Relaxation-Linearization Technique (RLT)}, a Second Order Cone Programming (SOCP) 
problem whose solution approximates the solution of the original nonconvex problem. 
Finally, the SOCP approximations are carefully used to obtain useful lower bounds at the 
top levels of the subset enumeration tree, where existing algorithms just pass through.

\subsection{Outline of the paper}

 We begin by showing, in Section \ref{S:Formulation}, the alternative formulations of the problem that permit 
to get rid of the order statistics involved in Problem \eqref{P:LTS}. More specifically, we show that the problem can be cast as 
a nonlinear mixed-integer program or a concave program, both particular cases of the best subset selection problem. 
Then, in Section \ref{S:linearization} we introduce approximate convexifications for products involving binary variables and nonlinear 
continuous terms when an upper bound for the continuous variables is not available. We analyse the validity of the relaxation obtained using 
an estimation of the bound on the continuous variables and their applicability in a branch and bound algorithm. 
In Section \ref{S:SDPbounds} we show how to obtain  an estimation of the bound on the continuous variables to be used for the 
approximate convexification using a known SDP reformulation of the concave maximization problem. 
Section \ref{S:SBB} describes the actual implementation of a branch and bound algorithm incorporating bounds from approximate convexifications.
In Section \ref{S:numerics} we present the results of a computational study showing the performance gains in a branch and bound algorithm due to 
SOCP bounds in a large set of problems. Finally, Section \ref{S:conclusion} concludes the paper with a discussion on further avenues 
for research in this subject.

\section{Formulation as a best subset problem}\label{S:Formulation}
Problem \eqref{P:LTS} can be written as a mixed integer nonlinear program using the fact 
that for arbitrary $r\in \RR^n$, if $r_{1:n}\leq r_{2:n}\leq \cdots \leq r_{n:n}$ denote its ordered components,
then
\begin{equation*}\label{E:minsum-concave}
\sum_{i=1}^h r_{i:n}=\min \left\{ \sum_{i=1}^n w_i r_{i} \left | w\in \{ 0,1\}^n, \sum_{i=1}^n w_i=h \right\}
 \right. ,
\end{equation*}
and
\begin{equation*}\label{E:L-scale-min-with-w}
\sum_{i=1}^h r(\beta)^2_{i:n}= \min_{w\in \CC_h} 
\sum_{i=1}^n w_i r_i(\beta)^2,
 \end{equation*}
where $\CC_h=\{ w\in \{ 0,1\}^n, \sum_{i=1}^n w_i=h\}$ is a representation of all the subsets of $\{1,...,n\}$ of size $h$.

Therefore our original Problem \eqref{P:LTS} is equivalent to the following  mixed-integer nonlinear programming problem:
\begin{equation}\label{E:nonlinear-MIP}
\begin{array}{r@{}l}
\min&\quad \sum\limits_{i=1}^{n} w_k r_k^2,\\
 s.t& \\
& r+X\beta=y,\\
& e^\top w=h,\\
w\in \{0,1 &\}^n, \beta\in\RR^d, r \in \RR^n,\\
\end{array}
\end{equation}
where $e$ is the $n\times 1$ vector of ones.

By splitting  the variables of Problem \eqref{E:nonlinear-MIP} we see that it is equivalent to 
\begin{equation}\label{P:marginal}
\min \{\ v(w) \ \mid   w\in \CC_h\} , 
\end{equation}
where $v(w)$ is the  value of the weighted LS problem
\begin{equation}\label{E:marginal}
 v(w)=\inf \left\{  \sum\limits_{k=1}^{n} w_k r_k^2:\ \beta \in\RR^d, r=y-X\beta \right\}
\end{equation}
obtained by minimizing over $\beta$ and $r$ for fixed $w\in\CC_h$ in \eqref{E:nonlinear-MIP}.
Hence, Problem \eqref{E:nonlinear-MIP} amounts to selecting the subset of $h$ observations with the least sum of squared
residuals.
The function $v$ defined in \eqref{E:marginal} is concave, therefore Problem \eqref{P:LTS} can be thought as 
a concave minimization problem. \citet{Giloni-Padberg2002} were the first to show this property, and used it to 
devise a local minimization procedure. \citet{Nguyen-Welsch2010} revisited this formulation and derived an SDP
formulation of the corresponding maximization problem. 
Unfortunately, the degeneracy of the feasible domain makes it difficult to apply concave minimization algorithms to problem \eqref{P:LTS}. 
In this paper we tackle the problem in the form given in \eqref{E:nonlinear-MIP}.

\subsection{The enumeration tree}\label{SS:tree}
\begin{figure}[th]
 \tikzstyle{hojad} =[circle, inner sep=1pt,draw]
\tikzstyle{hoja} =[circle, inner sep=1pt]
 \tikzstyle{activo} =[circle,draw,inner sep=1pt,fill=black!75]
\noindent\makebox[\textwidth]{
\begin{tikzpicture}[level distance=1cm,scale=0.75]


 \node [circle] (a1) at (1,0) {$1$}	[sibling distance=2cm]
		child {[sibling distance=.5cm]
			node[circle] (b12) {$2$} edge from parent[draw=none]
				child{node[hojad] (c123) {$3$}}
				child{node[hojad] (c124) {$4$}}
				child{node[hojad] (c125) {$5$}}
				child{node[hojad] (c126) {$6$}}
			}
		child {[sibling distance=.5cm]
			node[circle] (b13) {$3$} edge from parent[draw=none]
				child{node[hojad] (c134) {$4$}}
				child{node[hojad] (c135) {$5$}}
				child{node[hojad] (c136) {$6$}}
			}
		child[missing] {};    
\node[circle] (b14) at (2.5,-1cm) {$4$} [sibling distance=.5cm]
 				child{ node[hojad] (c145) {$5$} }
 				child{ node[hojad] (c146) {$6$} };
\node[circle] (b15) at (3.6,-1cm) {$5$} 
				child{node[hojad] (c156) {$6$}};
\node [circle] (a2) at (6.5,0) {$2$}	[sibling distance=1.5cm]
		child {[sibling distance=.5cm]
			node[circle] (b23) {{$3$}} edge from parent[draw=none]
				child{node[hojad] (c234) {$4$}}
				child{node[hojad] (c235) {$5$}}
				child{node[hojad] (c236) {$6$}}
			}
		child {[sibling distance=.5cm]
			node[circle] (b24) {$4$} edge from parent[draw=none]
				child{node[hojad] (c245) {$5$}}
				child{node[hojad] (c246) {$6$}	}
			}
			child[missing] {};
\node[circle] (b25) at (7.5,-1cm){$5$}
			child{node[hojad] (c256) {$6$}};

\node[circle] (a3) at (9,0) {$3$}  [sibling distance=0.9cm]	
		child {node[circle] (b34) {$4$} [sibling distance=0.5cm]
			child{node[hojad] (c345) {$5$}}
			child{node[hojad] (c346) {$6$}}
			}
		child {node[circle] (b35) {$5$}
			child{node[hojad] (c356) {$6$}}
			};
\node[circle] (a4) at (10.5,0) {$4$}	[sibling distance=0.8cm]
		child {node[circle] (b45) {$5$}
			child{node[hojad,inner sep=1pt] (c456) {$6$}}
			};

\draw  (a2) -- (b25.north);

\node[circle] (main) at (5,1.5) {$\emptyset$}
edge (a1)
edge (a2)
edge (a3)
edge (a4);

 \draw  (a1) -- (b12);
 \draw  (a1) -- (b13);
 \draw  (a1) -- (b14);
 \draw  (a1) -- (b15.north);
 \draw  (a2) -- (b23);
 \draw  (a2) -- (b24);

\node[] (k0) at (13,1.5) {$k=0$};
\node[] (k1) at (13,0) {$k=1$};
\node[] (k2) at (13,-1) {$k=2$};
\node[] (k2) at (13,-2) {$k=3$};
\end{tikzpicture}
  }


%
\caption{The BB tree for $n=6$ and $h=3$.}\label{F:BB Tree}
\end{figure}
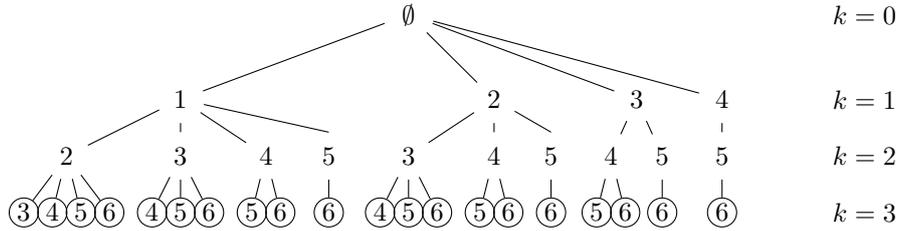

In Figure \ref{F:BB Tree} we depict the enumeration tree constructed by the branch and bound algorithm, in a small example with $h=3$ and $n=6$.

The circled nodes are the leaves; each leaf represents a subset of $3$ observations (an element of $\CC_h$), which is obtained adding 
recursively the parent of each node until the root $\emptyset$ is reached. For example, at the end of the second branch from right to 
left there are two leaves, the leaf at the right is associated to the subset of observations $6,4$ and $3$, and that at left to observations 
$5,4$ and $3$. 
In terms of the optimization variable $w$, they  are associated to the points $(0,0,1,1,0,1)^\top$ and $(0,0,1,1,1,0)^\top$ respectively.
Any node has associated two index sets $S_0$ and $S_1$ representing the variables fixed to $0$ and $1$ respectively, each of cardinality
$J_0$ and $J_1$. Using these two quantities we can compute the number of child nodes as $n-J_0-h+1$ and the number of leaves that can be 
reached from the \emph{l-th} child node as  $n-J_1-J_0-l \binom h-J_1+1$.

\subsection{The monotonocity lower bound}\label{SS:monotonicity}
The value of the function $v$ at a point $w\in [0,1]^n$ gives the least sum of \emph{weighted} squares of residuals with weights $w$.
In particular if $w\in \{0,1\}^n$, $w_i=1$ for $i\in J\subseteq\{1,...,n\}$ and $w_i=0$ for $i\notin J$, then $v(w)$ is the sum
of squares of the fit to the subset of observations $J$. The value of $v(w)$ is finite as long as the matrix $M(w)=X^\top D(w)X$ is invertible,
and in this case \citep{Agullo2001},

\begin{equation}\label{E:v(next)}
 v(w+ e_j)=v(w)+\frac{ r_j(w)^2}{1+ x_j^\top M(w)^{-1}x_j},
\end{equation}

where $D(w)$ is the diagonal matrix formed from the vector $w$, $r_j(w)$ is the $j$-th residual obtained from the weighted least squares 
problem \eqref{E:marginal} with weights $w$, and $e_j$ is the $j$-th euclidean basis vector. Formula \eqref{E:v(next)} gives the change in the sum of squared 
residuals by adding  observation $j$ to the fit $J$, provided that $J$ contains at least $d$ linearly independent $x_i$s.
Note  that $ v(w+ e_j)-v(w)\geq 0$, which means that the objective function is non-decreasing from one node to any of its 
children. Therefore, the RSS at one node is a lower bound for the objective function of its children. If the tree is 
examined using a depth-first search strategy, it is possible to keep a Cholesky factorization of $M(w)$ from which to perform 
rank-one updates in order to quickly compute the quantities $x_j^\top M(w)^{-1}x_j$. The shorthands of the monotonicity 
bound were already mentioned: Formula \eqref{E:v(next)} requires $M(w)=X^\top D(w)X$ to be invertible, thus it is not applicable
at the $d$ top levels of the tree, and it provides loose bounds when $J_1$ is far from $h$ (when there are few more than
$d$ observations).

\section{Pseudo-convexifications with approximate bounds}\label{S:linearization}
At each node of the tree  we need to (under) estimate the value of the problem:
\begin{equation}\label{P:bounding}
\begin{array}{ll}
\min& \sum\limits_{i=1}^{n} w_k r_k^2,\\
 s.t& \\
& r+X\beta=y,\\
& e^\top w=h,\\
&0\leq w\leq 1,\\
&\beta\in\RR^d, r \in \RR^n\\
\end{array}
\end{equation}
with the additional constraints 
\begin{equation}\label{E:node-constraints}
\begin{array}{c}
 w_k=1,\ k\in S_1\\
 w_k=0,\ k\in S_0\\
\end{array}
\end{equation}
for two index sets $S_1,S_0$ particular to each node. 
Suppose that  we had an upper bound $\Pi_k$ for each quadratic term $r_k^2$; then we 
can apply a linearization technique  \citep{Glover75,Adams-EtAl2004,Adams-Sherali90,Adams-Sherali93} to get rid of the product of the continuous 
term $r_k^2$ with the binary variable $w_k$. This is done by constructing a new continuous variable $u_k$ to replace each product $w_k r_k^2$, and 
adding a number of additional constraints in order to ensure  that the value of the variable $u_k$ equals the product $w_k r_k^2$ at any feasible 
point of the new problem. Applied to Problem \eqref{P:bounding}  this procedure yields to the following SOCP problem

\begin{equation}\label{P:SOCP1}
\begin{array}{l@{}lr}
\min\ & \  \sum_{k=1}^n u_k &\\
s.t\quad&   &  \\
\qquad \ r_k^2 & -\Pi_k(1-w_k) \leq  u_k, &  1\leq k \leq n\\
& w_k=1,&\qquad k\in S_1\\
& w_k=0,&\qquad k\in S_0\\
&\pto{e}{w} = h &\\
&r+X\beta =y,&\\
 u\in \RR_+^n,&\  r\in \RR^n, \beta\in \RR^d, w \in \{0,1\}^n.&
 \end{array}
\tag{P}
\end{equation}
\noindent

From the standard theory of \citet{Glover75}, if 
$\Pi_k\geq \bar{r}_k^2$, where $\bar{r}$ are the residuals at a solution to Problem \eqref{P:LTS}, then Problem \eqref{P:SOCP1} coincides with 
\eqref{P:LTS}  for any binary realization of $w$.  Unfortunately, as we show in the following section, the nonlinear coupling of the variables 
 restrains us from efficiently obtaining a guaranteed upper bound  for the residuals of Problem \eqref{P:LTS}. For this reason we introduce 
\emph{pseudo-convexifications}, which are instances of problem \eqref{P:SOCP1} for an approximate upper bound $(\Pi_k)_{1\leq k\leq n}$. 
We say that a solution to \eqref{P:SOCP1} is \emph{consistent} if $r_k^2<\Pi_k$ for any $k=1,...,n$ at the  optimal $r$.
Consistency of the solution to a pseudo-convexification is a necessary condition for being an actual convexification, but in general 
 it is not sufficient. As a consequence, using bounds obtained from a pseudo-convexification in a branch and bound algorithm can 
lead to pruning branches potentially yielding a new solution. Of course, this drawback can be avoided by giving very large values 
to $\Pi$. Nevertheless, larger values of $\Pi$ yields to weaker lower bounds when relaxing the binary constraint.
Therefore, we are led to find a compromise between a large bound ensuring the equivalence of Problems \eqref{P:LTS} and \eqref{P:SOCP1}, and 
a smaller bound promoting tight lower bounds for Problem \eqref{P:bounding}. 
In the next section we describe a strategy for obtaining bounds achieving that comprise.

\section{Obtaining Good Approximate bounds}\label{S:SDPbounds}

It is customary when linearizing polynomial programs to suppose that the optimization takes 
place on a bounded separable polytope, and that upper and lower bounds for each variable exists and can be computed by linear programming.
This does not hold in our case. In our problem the variable $w$ ranges over the unit cube in $\RR^n$, and all the other variables have a nonlinear 
dependency on $w$.
Indeed, for each $w$ there exists an unique feasible $\beta_w$, as seen from \eqref{E:marginal}; even more, $\beta_w$ is the unique solution
to the system $X^\top D(w)X \beta=X^\top D(w)y$. The residuals are linked to $\beta$, and therefore to $w$, by the linear constraint $r+X\beta =y$.
Therefore, an exact upper bound for $r_k^2$ can be obtained by solving the following auxiliary problem 
\begin{equation}\label{P:Ak}
\begin{array}{ll}
\max\limits_w &  r_k^2 \\
s.t& \\
& e^\top w=h\\
& r+X\beta =y,\\
&X^\top D(w) r=0\\
& 0\leq w \leq 1. \\
\end{array} \tag{$A_k$}
\end{equation}

Problem \eqref{P:Ak} is a maximization problem on $w$ only but, unlike Problem \eqref{P:marginal}, it is not a concave maximization problem.
Even it it could be done efficiently, solving $n$ problems would be cumbersome. For that reason 
we look for one single bound for all the residuals that can be efficiently computed.
A closely related problem is that of maximizing the weighted sum of squared residuals, 
which amounts to maximizing the function $v$ defined in \eqref{E:marginal},
\begin{equation}\label{P:RMTS}
\begin{array}{ll}
\max  & v(w)\\
s.t& \\
& e^\top w =q,\\
&0\leq w\leq 1,
\end{array} 
\end{equation}
for some $1\leq q \leq n$, This is a concave maximization problem, therefore a global minimum can be efficiently computed. In fact,  
\citet{Nguyen-Welsch2010} showed that Problem $\eqref{P:RMTS}$ can be cast as an SDP problem,
therefore it can be solved using standard widely-available software. If $d+1\leq q \leq n$ and we force 
the variables to be binary, the solution to Problem \eqref{P:RMTS} is the subset of $q$ observation with the largest sum 
of squared residuals. In any case, the value of the problem is monotone non-decreasing for $1\leq q \leq n$. Moreover, it has 
the property of going to $+\infty$ if any subset of the observations  is replaced by divergent ones.
After extensive numerical experiments we found that solving  $\eqref{P:RMTS}$ for $q=d/2$ yields to an effective approximation 
of the upper bound for the residuals.

\section{A branch and bound algorithm with SOCP bounds (S-BB)}\label{S:SBB}
In our implementation, the Narendra-Fukunaga tree described in  Subsection \ref{SS:tree} is examined using a depth-first search strategy.
We perform an in-level node ordering to take advantage of the unbalanced structure of the tree, as leftmost branches have 
many more children than those at the right. The innovations of our algorithm take place at the $d$ top levels of the tree;
at level $d+1$ and below the algorithm behaves like the BBA of \citet{Agullo2001}.

\subsection{Preliminaries}
We used the LS estimator as an initial solution; as a consequence, our algorithm is deterministic.
As indicated in Section \ref{S:SDPbounds}, we set $\Pi$ equal to the optimal value of Problem $\eqref{P:RMTS}$ with $q=d/2$.
Problem $\eqref{P:RMTS}$ is solved using the SDPT3 interior-point solver \citep{SDPT3}.  

\subsection{Lower bounds}
Problem \eqref{P:SOCP1} with the binary constraint relaxed to $0\leq w_k\leq 1$ is denoted as $\overline{\eqref{P:SOCP1}}$.  
The quadratic constraint $r_k^2-\Pi(1-w_k) \leq  u_k$ can be cast as a SOCP constraint, for that reason we call  
Problem $\overline{\eqref{P:SOCP1}}$ the SOCP relaxation. Problem $\overline{\eqref{P:SOCP1}}$ was solved using CPLEX 12.5.
At each node we first check the size of the subtree below each child node. 
Even nodes at top levels of the tree can have very few leaves below, in which case it is not worth spending time solving the SOCP relaxation. 
We launch the SOCP relaxation only for children with more than $10^6$ leaves.
If the solution  to Problem $\overline{\eqref{P:SOCP1}}$ results to be consistent, and the value of the problem is greater than the current upper 
bound, the branch is pruned. Otherwise, the residuals of the solution are still useful for ranking the children and performing the 
in-level ordering, by putting the observations with the largest residuals at the left, to promote subsequent pruning of large branches.

\subsubsection{Adjusting the bounds on $r$}
A look at problem \eqref{P:SOCP1} shows that the optimal values of $\Pi_k$ can be anticipated for $k\in S_0\cup S_1$.
For $k\in S_1$ the upper bound does not enter into play, we always have $u_k=w_k r_k^2$; on the contrary, for $k\in S_0$ 
we should always have $u_k=0$, therefore we set $\Pi_k=10\cdot \Pi$, where $\Pi$ is the upper bound obtained by solving Problem 
$\eqref{P:RMTS}$ and used for $k\notin S_0\cup S_1$. In practice this forces $u_k=0$, and does not spoil the conditioning of the 
SOCP problems as a huge number would do (in theory, we should set those $\Pi_k$ to $+\infty$).

\subsection{Local improvements}
Another innovation of our BB algorithm is the incorporation of a local search. Each time a leaf is examined, we apply
the \emph{concentration steps} of \citet{Rousseeuw-VanDriessen2006} to obtain an eventually better incumbent solution, and 
use it to update the global upper bound of the algorithm.

\section{Computational study}\label{S:numerics}
\begin{figure}[ht]
\begin{center}
\includegraphics[scale=0.7]{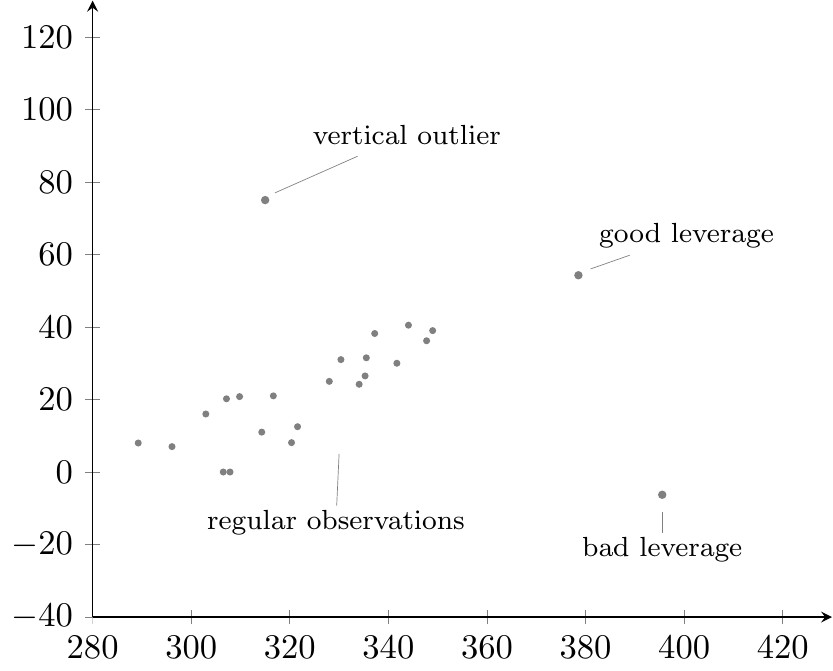}
\end{center}
\caption{A taxonomy of regression outliers.}\label{F:taxonomy}
 \end{figure}

Now we illustrate through a computational study the impact of incorporating the SOCP lower bound in a branch-and-bound algorithm.
In order to perform a systematic study, we generated synthetic datasets with sizes in a controlled range. 
It is known that the structure of the outliers, and not only their magnitudes, strongly affects the regression technique as 
well as the behaviour of the approximating algorithm. In Figure \ref{F:taxonomy}  we illustrate the taxonomy of 
linear regression outliers \citep{Rousseeuw-vanZomeren}. 
Outliers are called vertical if only the $y$ component (the response) is contaminated.
Vertical outliers are the more benign ones, and even some convex estimators, such as the $\ell_1$ estimator, can cope with them to some extent 
\citep{Giloni-Padberg2004}.
 Leverage points are points whose explicative variables are corrupted. In contrast to vertical outliers, 
leverage outliers can be very harmful. Excepting the case of the ``good leverage points'' 
illustrated in Figure  \ref{F:taxonomy}, which are in general not considered as outliers, bad leverage points 
are the most adversarial type of contamination.

For our study we generate synthetic data with outliers in the following way:
\begin{itemize}
 \item[-] We generate regular observations following model \eqref{E:linmod}, with $\delta$ standard Normal,
for different number of cases ($n$) and explicative variables ($d$).
\item[-] On each dataset we replace $10$ regular observations by bad leverage outliers.
\end{itemize}
The bad leverage outliers were obtained by shifting randomly selected observations in two different ways:
by a large, deterministic shift (high leverage outliers) and by adding a random term drawn form a Laplace distribution 
(heavy tail outliers).

For each combination $n/d/$type-of-contamination we drew $25$ datasets as described above and 
measured the total time spent by the S-BB algorithm described in Section \ref{S:SBB} and by the BBA \citep{Agullo2001} in 
computing the LTSE with a breakdown point of $50\%$ ($h=\pentiere{n/2}+\pentiere{(d+1)/2}$).

{
All computations were done in MATLAB version R2008b on a 64-bit Linux machine, with 8 cores and 6 GB RAM. 
For the SOCP relaxations we used IBM ILOG CPLEX Studio v. 12.5 via its MATLAB interface.}

The computing  times in tens of seconds, averaged over the $25$ repetitions, are shown in Tables \ref{T:Numtest} and \ref{T:NumtestHT} 
respectively. The impact of the SOCP bounds is, as expected, more important as the number of explicative variables increases, and 
more pronounced for larger $n$. The reduction in computing time exceeds the $20\%$ for $n=40$ and $d$ greater than $15$.
The accuracy of the solution is largely preserved; in Table \ref{T:success} we show the rate of success, which is 
larger than $99\%$  in all but one of the cases. The computing times in Tables \ref{T:Numtest} and \ref{T:NumtestHT} 
 marked with a dagger are averages excluding the run that did not gave the exact solution.

Further reductions in computing time are possible by relaxing the optimality goal. In this direction, we can 
mention that the fraction of the SOCP relaxations resulting in inconsistent solutions is not negligeable; 
using those solutions to derive bounds could result in a great performance improvement.  
Another way to do the same thing is by decreasing the parameter $q$ used to obtain the approximate bound $\Pi$.
However, the goal of this work was to improve the computation of the LTSE with the least possible 
lost in accuracy, and it was achieved.

\begin{table}[bh]
\begin{center}
\begin{tabular}{@{\extracolsep{-2pt}}llccccccc}
               &       &\multicolumn{7}{c}{\it d} \\ 
Alg.           & $n$               & $12$        & $13$     &  $14$ & $15$ &$16$&$17$&$18$   \\ \hline
\multirow{3}{*}{BBA}
   & $30$ &      0.51 &   0.35  &  0.51  &  0.26 &  0.36  &  0.14   & 0.18\\
   & $35$ &       5.19  &   4.79   &  7.99  &  5.68  &  8.72  &  5.02  &  7.29   \\
   & $40$ &    22.15 &  23.68 &  40.85 &  33.59 &  56.63 &  38.99 &  59.27   \\ \hline
\multirow{3}{*}{S-BB}
   & $30$ &      0.49  &  0.34  &  $0.51^\dagger$  &  0.26 &   0.35 &   0.14 &   0.18   \\
   & $35$ &      4.60   &  4.21  &  6.80  & $5.05^\dagger$  & $7.65^\dagger$  &  4.45 & $6.41^\dagger$   \\
   & $40$ &    22.53 &  20.89 &  34.90 &  26.91 &  45.44 &  29.87 &  46.49
\end{tabular}
\caption{Computing times (in tens of seconds) for high leverage outliers}\label{T:Numtest}
\end{center}
\end{table}
\begin{table}
\begin{center}
\begin{tabular}{@{\extracolsep{-2pt}}llccccccc}
               &       &\multicolumn{7}{c}{\it d} \\ 
Alg.           & $n$        & $12$        & $13$     &  $14$ & $15$ &$16$&$17$&$18$   \\ \hline
\multirow{3}{*}{BBA}
   & $30$ &   0.52  &  0.35  &  0.52  &  0.26  &  0.35 &   0.14 &   0.19   \\
   & $35$ &   5.26  &  4.81   & 8.07   & 5.67  &  8.69  &  5.03  &  7.29       \\
   & $40$ &  21.41 &  23.90 &  41.41 &  34.82 &  56.57  & 39.71 &  59.47   \\ \hline
\multirow{3}{*}{S-BB}
   & $30$ &   0.50  &  0.34 &   0.52  &  $0.26^\dagger$ &   0.35  &  0.14 &   0.19 \\
   & $35$ &   4.52 &  4.21 &   6.86  &  4.83  & 7.71  & $4.48^\dagger$   &  6.53     \\
   & $40$ &   22.11 &  21.43 &  34.70 &  27.46 &  $44.54^\dagger$ &  31.44 &  46.72  
\end{tabular}
\caption{Computing times (in tens of seconds) for heavy tail outliers}\label{T:NumtestHT}
\end{center}
\end{table}

\begin{table}
\begin{center}
\begin{tabular}{@{\extracolsep{-2pt}}ccccr}
                 &\multicolumn{3}{c}{\it n} \\ 
	$\%$		 &  $30$  & $35$ 	&  $40$   & average      \\ \hline
{high leverage outliers} & $99.55$   & $98.66$ & $100$   &   $99.40$              \\
 {heavy tail outliers}   & $99.55$ &   $99.55$ & $99.55$   &  $99.55$  \\
{average}                &$99.55$  & $99.10$ & $ 99.77$  & $99.47$\\
\end{tabular}
\caption{Success rates of the S-BB algorithm}\label{T:success}
\end{center}
\end{table}

\section{Conclusions and perspectives}\label{S:conclusion}

We have presented an approximate convex relaxation for the LTS problem. Its incorporation in a branch-and-bound algorithm
yields to significant savings in computing time at the price of a negligeable accuracy lost.
Our standpoint is that of improving the computation of an estimator with well-studied statistical properties. 
Now then, with the understanding of the underlying problem gained with this study one could propose alternative
techniques yielding to robust estimators defined with an eye on computability. Concretely, think of 
the estimator defined as the solution to Problem \eqref{P:SOCP1}, with $\Pi$ obtained by solving 
\eqref{P:RMTS} for some $1\leq q \leq n$. If the binary constraint is relaxed, such an estimator 
is defined by two convex optimization problems, with a range of applicability in the thousands 
of observations; if the binary constraint is kept, we dispose of SOCP convex relaxations at any node 
of the BB tree, without wasting time with inconsistent relaxations. The study of the statistical and 
computational aspects of that proposal will be the subject of a future work.

\bibliographystyle{elsarticle-harv}
\bibliography{articles,biblio}

\begin{thebibliography}{28}
\expandafter\ifx\csname natexlab\endcsname\relax\def\natexlab#1{#1}\fi
\expandafter\ifx\csname url\endcsname\relax
  \def\url#1{\texttt{#1}}\fi
\expandafter\ifx\csname urlprefix\endcsname\relax\def\urlprefix{URL }\fi

\bibitem[{Adams et~al.(2004)Adams, Forrester, and Glover}]{Adams-EtAl2004}
Adams, W.~P., Forrester, R.~J., Glover, F.~W., 2004. Comparisons and
  enhancement strategies for linearizing mixed 0-1 quadratic programs. Discrete
  optimization 1~(2), 99--120.

\bibitem[{Adams and Sherali(1990)}]{Adams-Sherali90}
Adams, W.~P., Sherali, H.~D., 1990. Linearization strategies for a class of
  zero-one mixed integer programming problems. Operations Research 38~(2),
  217--226.

\bibitem[{Adams and Sherali(1993)}]{Adams-Sherali93}
Adams, W.~P., Sherali, H.~D., 1993. Mixed-integer bilinear programming
  problems. Mathematical Programming 59~(3), 279--305.

\bibitem[{Agull{\'o}(2001)}]{Agullo2001}
Agull{\'o}, J., 2001. New algorithms for computing the least trimmed squares
  regression estimator. Computational Statistics \& Data Analysis 36~(4),
  425--439.

\bibitem[{Bernholt(2005)}]{Thorsten}
Bernholt, T., 2005. Computing the least median of squares estimator in time
  {O}($n^d$). In: Gervasi, O., Gavrilova, M.~L., Kumar, V., Lagana, A., Lee,
  H.~P., Mun, Y., Taniar, D., Tan, C. (Eds.), Computational Science and Its
  Applications -- {ICCSA} 2005. Vol. 3480 of Lecture Notes in Computer Science.
  Springer Berlin Heidelberg, pp. 697--706.

\bibitem[{Bertsimas and Mazumder(2014)}]{Bertsimas}
Bertsimas, D., Mazumder, R., 12 2014. Least quantile regression via modern
  optimization. Ann. Statist. 42~(6), 2494--2525.

\bibitem[{Chen({2003})}]{CHEN2003}
Chen, X., {AUG} {2003}. {An improved branch and bound algorithm for feature
  selection}. {Pattern Recognition Letters} {24}~({12}), {1925--1933}.

\bibitem[{Donoho and Huber(1983)}]{Donoho-Huber83}
Donoho, D., Huber, P.~J., 1983. The notion of breakdown point. In: A
  {F}estschrift for {E}rich {L}. {L}ehmann. Wadsworth Statist./Probab. Ser.
  Wadsworth, pp. 157--184.

\bibitem[{Erickson et~al.(2006)Erickson, Har-Peled, and Mount}]{Mount2006}
Erickson, J., Har-Peled, S., Mount, D.~M., 2006. On the least median square
  problem. Discrete Comput. Geom. 36~(4), 593--607.

\bibitem[{Giloni and Padberg(2002)}]{Giloni-Padberg2002}
Giloni, A., Padberg, M., 2002. Least trimmed squares regression, least median
  squares regression, and mathematical programming. Mathematical and Computer
  Modelling 35~(9-10), 1043 --1060.

\bibitem[{Giloni and Padberg(2004)}]{Giloni-Padberg2004}
Giloni, A., Padberg, M., 2004. The finite sample breakdown point of
  $\ell_1$-regression. SIAM Journal on Optimization 14~(4), 1028--1042.

\bibitem[{Glover(1975)}]{Glover75}
Glover, F., 1975. Improved linear integer programming formulations of nonlinear
  integer problems. Management Science 22~(4), 455--460.

\bibitem[{Hofmann et~al.(2010)Hofmann, Gatu, and Kontoghiorghes}]{Gatu2010}
Hofmann, M., Gatu, C., Kontoghiorghes, E.~J., 2010. An exact least trimmed
  squares algorithm for a range of coverage values. Journal of Computational
  and Graphical Statistics 19~(1), 191--204.

\bibitem[{Jalali-Heravi and Konouz(2002)}]{krafft}
Jalali-Heravi, M., Konouz, E., 2002. Use of quantitative structure-property
  relationships in predicting the {K}rafft point of anionic surfactants.
  Electron. J. Mol. Des. 1, 410--417.

\bibitem[{Maronna et~al.(2006)Maronna, Martin, and
  Yohai}]{Maronna-Martin-Yohai2006}
Maronna, R.~A., Martin, R.~D., Yohai, V.~J., 2006. Robust statistics. Wiley
  Series in Probability and Statistics. John Wiley \& Sons.

\bibitem[{Mount et~al.(2014)Mount, Netanyahu, Piatko, Silverman, and
  Wu}]{Mount2014}
Mount, D.~M., Netanyahu, N.~S., Piatko, C.~D., Silverman, R., Wu, A.~Y., 2014.
  On the least trimmed squares estimator. Algorithmica 69~(1), 148--183.

\bibitem[{Mount et~al.(2007)Mount, Netanyahu, Romanik, Silverman, and
  Wu}]{Mount2007}
Mount, D.~M., Netanyahu, N.~S., Romanik, K., Silverman, R., Wu, A.~Y., 2007. A
  practical approximation algorithm for the {LMS} line estimator. Comput.
  Statist. Data Anal. 51~(5), 2461--2486.

\bibitem[{Narendra and Fukunaga(1977)}]{Narendra-Fukunaga}
Narendra, P., Fukunaga, K., 1977. A branch and bound algorithm for feature
  subset selection. IEEE Transactions on Computers C-26~(9), 917 --922.

\bibitem[{Nguyen and Welsch(2010)}]{Nguyen-Welsch2010}
Nguyen, T.~D., Welsch, R., 2010. Outlier detection and least trimmed squares
  approximation using semi-definite programming. Comput. Statist. Data Anal.
  54~(12), 3212--3226.

\bibitem[{Rousseeuw and {Van Driessen}(2006)}]{Rousseeuw-VanDriessen2006}
Rousseeuw, P., {Van Driessen}, K., 2006. Computing {LTS} regression for large
  data sets. Data Mining and Knowledge Discovery~(12), 29--45.

\bibitem[{Rousseeuw and Leroy(1987)}]{Rousseeuw-Leroy87}
Rousseeuw, P.~J., Leroy, A.~M., 1987. Robust regression and outlier detection.
  John Wiley \& Sons.

\bibitem[{Rousseeuw and van Zomeren(1990)}]{Rousseeuw-vanZomeren}
Rousseeuw, P.~J., van Zomeren, B.~C., 1990. Unmasking multivariate outliers and
  leverage points. Journal of the American Statistical Association 85~(411),
  633--639.

\bibitem[{Somol et~al.({2004})Somol, Pudil, and Kittler}]{Somol2004}
Somol, P., Pudil, P., Kittler, J., {JUL} {2004}. {Fast branch \& bound
  algorithms for optimal feature selection}. {IEEE Transactions On Pattern
  Analysis And Machine Intelligence} {26}~({7}), {900--912}.

\bibitem[{Steele and Steiger(1986)}]{Steele-Steiger1986}
Steele, J.~M., Steiger, W.~L., May 1986. Algorithms and complexity for least
  median of squares regression. Discrete Appl. Math. 14~(1), 93--100.

\bibitem[{Stromberg(1993)}]{Stromberg1993}
Stromberg, A.~J., Nov. 1993. Computing the exact least median of squares
  estimate and stability diagnostics in multiple linear regression. SIAM J.
  Sci. Comput. 14~(6), 1289--1299.

\bibitem[{Torti et~al.(2012)Torti, Perrotta, Atkinson, and
  Riani}]{TortiEtAl2012}
Torti, F., Perrotta, D., Atkinson, A.~C., Riani, M., 2012. Benchmark testing of
  algorithms for very robust regression: {FS}, {LMS} and {LTS}. Computational
  Statistics \& Data Analysis 56~(8), 2501--2512.

\bibitem[{T\"{u}t\"{u}nc\"{u} et~al.(2003)T\"{u}t\"{u}nc\"{u}, Toh, and
  Todd}]{SDPT3}
T\"{u}t\"{u}nc\"{u}, R.~H., Toh, K.~C., Todd, M.~J., 2003. Solving
  semidefinite-quadratic-linear programs using {SDPT3}. Mathematical
  Programming 95, 189--217.

\bibitem[{Yu and Yuan(1993)}]{YUAN1993}
Yu, B., Yuan, B., 1993. A more efficient branch-and-bound algorithm for
  feature-selection. Pattern Recognition 26~(6), 883--889.

\end{thebibliography}
\end{document}